\documentclass[12pt,oneside]{article}

\usepackage{amsfonts}
\usepackage{amscd}
\usepackage{amsthm}
\usepackage{amssymb}
\usepackage{amsmath}
\usepackage{epsfig}

\usepackage{graphicx}
\usepackage{xypic}
\input xypic
\input xy
\xyoption{all}

\newcommand{\mR}{\mathbb R}
\newcommand{\mZ}{\mathbb Z}
\newcommand{\mI}{\mathbb I}
\newcommand{\mH}{\mathbb H}
\newcommand{\mN}{\mathbb N}

\newcommand{\tX}{\tilde X}

\newcommand{\G}{\mathcal G}
\newcommand{\B}{\mathcal B}

\newcommand{\C}{\mathcal C}
\newcommand{\U}{\mathcal U}

\newcommand{\g}{\gamma}

\newcommand{\dlim}{\underrightarrow{\lim}}
\newcommand{\pinf}{\partial ^\infty}

\title{Diagram rigidity for geometric amalgamations of free groups.}

\author{Jean-Fran\c cois Lafont}

\theoremstyle{definition}
\newtheorem{Def}{Definition}[section]

\theoremstyle{proposition}
\newtheorem{Lem}{Lemma}[section]
\newtheorem{Prop}{Proposition}[section]

\theoremstyle{plain}

\newtheorem{Thm}{Theorem}[section]
\newtheorem{Cor}{Corollary}[section]

\theoremstyle{remark}

\newtheorem*{Prf}{Proof}

\begin{document}

\maketitle

\begin{abstract}
In this paper, we show that a class of 2-dimensional locally CAT(-1) spaces is
topologically rigid: isomorphism of the fundamental groups is equivalent to the
spaces being homeomorphic.  An immediate application of this result is a diagram
rigidity theorem for certain amalgamations of free groups.  The direct limit of
two such amalgamations are isomorphic if and only if there is an isomorphism
between the respective diagrams.  
\end{abstract}

\section{Introduction.}

In this paper, we introduce the notion of a {\it geometric amalgamation} of
free groups.  This is a class of diagrams of groups, with the property that
they are {\it rigid} in the following sense:

\begin{Thm}[Diagram rigidity]
Let $\G_1$, $\G_2$ be a pair of geometric amalgamations of free groups.  Then
$\dlim \G_1$ is isomorphic to $\dlim \G_2$ if and only if $\G_1$ is isomorphic
to $\G_2$ (as diagrams of groups).
\end{Thm}

In order to prove this theorem, we will first translate the question to a more
topological setting.  Associated to any geometric amalgamation of free groups,
there is a canonically defined topological space, which we call a {\it simple,
thick, 2-dimensional hyperbolic P-manifold}.  The associated space will have
fundamental group isomorphic to the direct limit of the geometric amalgamation,
and has the property that the diagram can be ``read off'' from the topology
of the space.  The first theorem will then be a consequence of the following
purely topological result:

\begin{Thm}[Topological rigidity]
Let $X_1$, $X_2$ be a pair of simple, thick, 2-dimensional
hyperbolic P-manifolds, and assume that
$\phi:\pi_1(X_1)\rightarrow \pi_1(X_2)$ is an isomorphism.  Then
there exists a homeomorphism $\Phi:X_1\rightarrow X_2$ that
induces $\phi$ on the level of the fundamental groups.
\end{Thm}

Immediate applications of this second theorem are the following:

\begin{Cor}[Nielson realization]
Let $X$ be a simple, thick, 2-dimensional hyperbolic P-manifold.  Then the
canonical map from $Homeo(X)$ to $Out(\pi_1(X))$ is surjective.
\end{Cor}

\begin{Cor}[co-Hopf property]
Let $\G$ be a geometric amalgamation of free groups.  Then $\dlim \G$ is a
co-Hopfian group.
\end{Cor}

We note that the second corollary can also be obtained from work of Sela 
\cite{S}, who proved that a non-elementary torsion-free $\delta$-hyperbolic
groups is co-Hopfian if and only if it is freely indecomposable.

The proof of Theorem 1.2 relies on a topological characterization of certain 
points in the boundary at infinity of the universal cover of a simple, thick,
2-dimensional hyperbolic P-manifold (these spaces are CAT(-1)).

\section{Some preliminaries.}

In this section, we briefly define the various notions that are relevant to this paper,
and recall some basic facts that will be used in the proofs.

\begin{Def}
For our purposes, a {\it diagram of groups} $\G$ will consist of:
\begin{itemize}
\item a finite connected directed graph, which we will also denote by $\G$, with vertex set
$V(\G)$ and directed edges $E(\G)$,
\item an assignement of a group $G_v$ to each vertex $v\in V(\G)$
\item an assignement of a homomorphism $\phi_e:G_{e_-}\rightarrow G_{e_+}$ to each
directed edge $e\in E(\G)$, where $e_-$ and $e_+$ denote the initial and terminal
endpoints of the directed edge $e$.
\end{itemize}
We will denote by $\dlim \G$ the direct limit of the diagram $\G$. 

\end{Def}

Observe that the above definition differs superficially from the
notion of a graph of groups used in Bass-Serre theory (see for
instance Serre \cite{Se}).  The interested reader can easily
translate the above definition into the language of Bass-Serre
theory.  We now refine the above definition to the diagrams we are
really interested in:

\begin{Def}
We will say that a diagram of groups $\G$ is a {\it geometric amalgamation of free groups}
provided it satisfies the following properties:
\begin{itemize}
\item the vertex set $V(\G)$ can be partitioned into $V_0(\G)$ and
$V_1(\G)$, and every directed edge $e\in E(\G)$ has $e_-\in
V_0(\G)$, $e_+\in V_1(G)$.  Furthermore, each vertex in $V_0(\G)$
has degree at least three. \item the group associated to every
$v\in V_0(\G)$ is isomorphic to $\mZ$, and the group associated to
every $w\in V_1(\G)$ is a free group with rank $\geq 2$. \item
associated to every $w\in V_1(\G)$, there is a compact surface
$M_w$ whose fundamental group is $G_w$. \item each edge morphism
$\phi_e$ maps $G_{e_-}\cong \mZ$ isomorphically onto a group
conjugate to the fundamental group of a boundary component in
$G_{e_+}\cong \pi_1(M_{e_+})$. \item for every vertex $w\in
V_1(G)$, and every conjugacy class of $\mZ$-subgroups of $G_w\cong
\pi_1(M_w)$ corresponding to a boundary component of $M_w$, there
is precisely one edge $e\in E(\G)$ with $e_+=w$, and
$\phi_e(G_{e_-})$ lying within the conjugacy class.
\end{itemize}
\end{Def}

More concisely, we can think of a geometric amalgamation of free groups as being a diagram
of groups consisting of two rows:
$$
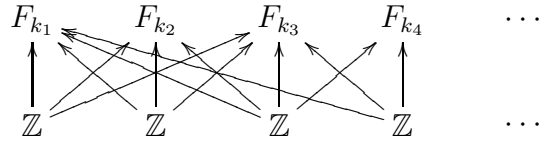
\xymatrix{ F_{k_1} & F_{k_2} & F_{k_3} & F_{k_4} &  \cdots\\
\mZ \ar[u] \ar[ur] \ar[urr] & \mZ \ar[ul] \ar[u] \ar[ur] & \mZ \ar[u] \ar[ull] \ar[ul]
\ar[ur] & \mZ \ar[u] \ar[ul] \ar[ulll] & \cdots }
$$
where each $\mZ$ group in the bottom row injects into at least
three of the free groups in the top row, and has image lying in a
``boundary subgroup'' of the free groups in the top row.  In
addition, each ``boundary subgroup'' in the top row lies in the
image of precisely one $\mZ$ from the bottom row (upto conjugacy).

Now the main motivation behind our terminology lies in that the direct limit of a
geometric amalgamation of free groups naturally corresponds to the fundamental group
of an associated topological space.  We now proceed to define these spaces.

\begin{Def}
We say that a compact geodesic metric space $X$ is a simple, thick, 2-dimensional hyperbolic
P-manifold, provided that there exists a closed subset $Y\subset X$ with the property that:
\begin{itemize}
\item each connected component of $Y$ is homeomorphic to $S^1$, and forms a totally geodesic
subspace of $X$.
\item the closure of each connected component of $X-Y$ is homeomorphic to a compact
surface with boundary, and the homeomorphism takes the component of $X-Y$
to the interior of the surface with boundary; the closure of such a
component will be called a {\it chamber}.
\item there exists a hyperbolic Riemannian metric on each chamber which coincides with the
original metric.
\item each connected component of $Y$ lies in the closure of at least three distinct
chambers.
\end{itemize}
We will call the subset $Y$ the {\it branching locus}, and will call the connected
components of $Y$ {\it branching geodesics}.
\end{Def}

Note that simple, thick, 2-dimensional hyperbolic P-manifolds are
locally CAT(-1) (see \cite{bh}), and hence their universal covers
are CAT(-1) spaces.  In particular, this implies that their
fundamental groups are $\delta$-hyperbolic groups, and that an
abstract isomorphism between the fundamental groups of two such
spaces naturally induces a quasi-isometry between their universal
covers.

To make precise the correspondance between the previous two definitions, we show the
following:

\begin{Lem}
A group $G$ is the fundamental group of a simple, thick, 2-dimensional hyperbolic P-manifold
if and only if it is the direct limit of a geometric amalgamation of free groups.
\end{Lem}

\begin{Prf}
If $G$ is the fundamental group of a simple, thick, 2-dimensional
hyperbolic P-manifold $X$, we consider an open cover $\{\U_i\}$ of
$X$ by open sets which consist of $\epsilon$-neighborhoods of the
chambers of $X$, where $\epsilon$ is chosen to be small enough. An
immediate application of the general form of the Siefert-Van
Kampen Theorem (see Chapter 2, Section 7 in May \cite{ma}) that
the fundamental group of $X$ is isomorphic to the direct limit of
a diagram of groups obtained from the covering $\{\U_i\}$.
Furthermore, the diagram has vertex groups isomorphic to the
fundamental groups of the various intersections of open sets in
the covering $\{\U_i\}$, with edge morphisms induced by
inclusions.  It is immediate that for the covering we've defined,
the resulting diagram is a geometric amalgamation of free groups.
Furthermore, the diagram is uniquely defined by the space $X$.

Conversely, assume that $G$ is the direct limit of a geometric
amalgamation of free groups, denoted by $\G$.  Corresponding to
the diagram $\G$, we can associate a diagram of topological spaces
by associating to each vertex in $v\in V_0$ an $S^1$, and to each
vertex $w\in V_1$ the corresponding compact surface with boundary
$M_w$.  Note that to each edge, there corresponds a homeomorphism
from one of the $S^1$ (corresponding to the initial vertex of the
edge) to a boundary component of one of the $M_w$ (corresponding
to the terminal vertex of the edge).  Now consider the direct
limit of this diagram of spaces in the category of topological
spaces.  It is immediate that this direct limit is a simple,
thick, 2-dimensional P-manifold $X$, and, by the discussion in the
previous paragraph, that $\pi_1(X)\cong G$.

To conclude, we need to show that $X$ supports a hyperbolic
metric.  To see this, we make each $S^1$ isometric to the unit
circle in $\mR^2$, and make each $M_w$ isometric to a compact
hyperbolic surface with all boundary components totally geodesic
of length $2\pi$.  We further require the homeomorphisms from the
$S^1$ to the boundary components of the $M_w$ to be isometries.
This immediately yields a hyperbolic metric on the P-manifold $X$.

Finally, we point out that the space $X$ constructed above is
unique upto homeomorphism.  This follows from the fact that if $M$
is a compact orientable surface with boundary, and $\phi:\partial
M\rightarrow \partial M$ is an orientation preserving
self-homeomorphism, then there is a self-homeomorphism $\hat
\phi:M\rightarrow M$ that induces $\phi$ when restricted to the
boundary.  In particular, the choice of homeomorphisms used to identify
the various $S^1$ (corresponding to the $\mZ$ groups) with the boundary
components of the $M_w$ (corresponding to the free groups) does not 
influence the topology of the resulting space. 
\end{Prf}

Finally, to conclude this section, we reduce the proof of Theorem
1 to that of Theorem 2:

\begin{Prf}[Diagram Rigidity]
Let us start with a pair $\G_1$, $\G_2$ of geometric amalgamations
of free groups, and assume that $\dlim \G_1$ is isomorphic to
$\dlim \G_2$.  From the previous Lemma, we can associate to each
$\G_i$ a simple, thick, 2-dimensional hyperbolic P-manifold $X_i$,
with the property that $\pi_1(X_i)\cong \dlim \G_i$.  In
particular, the isomorphism between direct limits yields an
isomorphism $\phi:\pi_1(X_1)\rightarrow \pi_1(X_2)$. From the
topological rigidity result, there is a homeomorphism
$\Phi:X_1\rightarrow X_2$ which induces $\phi$ on the level of
fundamental groups.

Now note that $\Phi$, being a homeomorphism, must map the
branching locus in $X_1$ homeomorphically to the branching locus
in $X_2$, and hence maps the chambers of $X_1$ homeomorphically to
the chambers of $X_2$.  Furthermore, the map $\Phi$ induces a
bijection between the chambers in $X_1$ and those in $X_2$.  But
by the uniqueness portion of the Lemma above, this implies that
the diagrams $G_1$ and $G_2$ are isomorphic, concluding our proof.
\end{Prf}


\section{Topological rigidity.}

In this section, we provide a proof of Theorem 1.2.
Let us start by fixing some notation. $X$ will always denote a
simple, thick, 2-dimensional hyperbolic P-manifolds, $\tX$ the
universal cover of $X$, and $\pinf \tX$ the boundary at
infinity of $\tX$.  We will let $\Gamma$ denote the
fundamental group of $X$.  $\G$ will denote the collection of
geodesics in $\tX$, and $\B\G\subset \G$ will denote the
collection of lifts of branching geodesics in $\tX$.  We will
let $\C$ denote the collection of lifts of chambers in $\tX$.
Finally, we will let $\pinf \B\G \subset \pinf \tX$ be the
collection of points of the form $\g (\pm \infty)$ where $\g \in
\B\G$.  In the portions of this section where we deal with a pair
of simple, thick, 2-dimensional hyperbolic P-manifolds, we will
append subscripts to keep track of which of the two spaces
we are referring to.

The first step in our argument consists of identifying the
separation properties of individual points in $\pinf \tX$.  We
start with an easy:

\begin{Lem}
The boundary at infinity $\pinf \tX$ is path-connected.
\end{Lem}

\begin{Prf}
To see this, let $p_1,p_2\in \pinf \tX$ be an arbitrary pair of
points.  Let $\gamma$ be a geodesic in $\tX$ satisfying
$\gamma(-\infty)=p_1$, $\gamma(\infty)=p_2$.  The simplicity and
thickness hypotheses on $X$ ensure that there exists an
isometrically embedded $f:\mH^2\hookrightarrow \tX$ with the
property that $\gamma\subset \mH^2$.  In particular, we see that
the pair of points $p_1,p_2$ lie on an embedded $S^1=\pinf \mH^2
\hookrightarrow \pinf \tX$.  Hence there exists a path in $\pinf
\tX$ joining $p_1$ to $p_2$, concluding the proof of the Lemma.
\end{Prf}

Let us recall some basic definitions.  A subset $S$ in a
topological space $X$ is said to {\it locally separate} provided
there exists a neighborhood $N$ of $S$ with the property that
$N-S$ is disconnected.  If there exists a neighborhood $N$ such
that $N-S$ consists of $\geq m$ connected components, we say that
$S$ locally separates into $\geq m$ components.  We say that $S$
locally separates into $M$ components provided $M$ is the supremum
of the integers $m$ with the property that $S$ locally separates
into $\geq m$ components.  We will be interested in the case where
$S$ consists either of a single point, or of a pair of points.

In the case where $S=\{x\}$, and $S$ locally separates the space
$X$, we will say that $x$ is a {\it local cutpoint} of $X$.  If in
addition $X-\{x\}$ is disconnected,
will say that the point $x$ is a {\it global cutpoint} of $X$. 



We observe that the previous Lemma in particular implies that
$\pinf \tX$ is connected.  From the work of Bowditch \cite{bo}
and Swarup \cite{sw}, this immediately yields:

\begin{Cor}
The boundary at infinity $\pinf \tX_i$ does not contain any global
cutpoints.
\end{Cor}

We note that the previous Corollary tells us that the {\it global}
separation properties of individual points in $\pinf \tX_i$ are
uninteresting.  On the other hand, the {\it local} separation
properties of points in $\pinf \tX_i$ are quite interesting.  Our
next step is to consider a notion which is slightly weaker than
``local separation into $\geq 3$ components''.

We now collect some basic facts concerning separation and
connectedness properties in simple, thick, hyperbolic P-manifolds.
The proofs of the following four lemmas can be found in \cite{la1}
in the 3-dimensional setting, but the arguments given there extend
verbatim to the 2-dimensional setting.

\begin{Lem}[Lemma 2.1 in \cite{la1}]
Let $\g \in \B\G$ be a branching geodesic in $\tX$, and let
$C_1,C_2\in \C$ be two lifts of chambers which are both incident
to $\g$. Then $C_1-\g $ and $C_2-\g$ lie in different connected
components of $\tX - \g$.
\end{Lem}

\begin{Lem}[Lemma 2.2 in \cite{la1}]
Let $C\in \C$ be a lift of a chamber, and let $\g_1,\g_2\in \B\G$
be two branching geodesics which are both incident to $C$. Then
$\g_1$ and $\g_2$ lie in different connected components of $\tX -
Int(C)$.
\end{Lem}

\begin{Lem}[Lemma 2.3 in \cite{la1}]
Let $\{\g(\pm\infty)\}$ be the pair of points in $\pinf \tX$
corresponding to some $\g\in \B\G$, and let $\pinf C_1$, $\pinf
C_2$ be the boundaries at infinity of two lifts of chambers
$C_1,C_2\in \C$ which are both incident to $\g$. Then $\pinf
C_1-\{\g(\pm\infty)\}$ and $\pinf C_1-\{\g(\pm\infty)\}$ lie in
different connected components of $\pinf \tX-\{\g(\pm\infty)\}$.
\end{Lem}

\begin{Lem}[Lemma 2.4 in \cite{la1}]
Let $\pinf C$ be the boundary at infinity corresponding to a
connected lift of a chamber $C\in \C$, and let
$\{\g_1(\pm\infty)\}$, $\{\g_2(\pm\infty)\}$ be the boundary at
infinity of two branching geodesics $\g_1,\g_2\in \B\G$ which are
both incident to $C$. Then $\{\g_1(\pm\infty)\}$ and
$\{\g_2(\pm\infty)\}$ lie in different connected components of
$\pinf \tX - (\pinf C - \cup \pinf \eta_i)$, where the union is
over all $\eta_i\in \B\G$ which are boundary components of $C$.
\end{Lem}

We define the {\it tripod} to be the space $T$ obtained by taking
the cone of a 3-point set.  The cone point will be denoted by
$*\in T$.  We say that a point $x$ is a {\it branching point} in a
topological space $X$ provided there exists an injective map
$f:T\rightarrow X$ satisfying $f(*)=x$. It is easy to see that if
a point $x$ in a geodesic space $X$ locally separates into $\geq
3$ components, then the point $x$ is a branching point in $X$. The
notion of branching was introduced by the author in \cite{la1},
\cite{la2} in order to study the local topology of the boundary at
infinity of simple, thick hyperbolic P-manifolds.  The proof of
the following Proposition closely parallels the arguments given in the paper
\cite{la1}:

\begin{Prop}
For a point $p\in \pinf \tX$, we have that $p$ is branching if and
only if $p=\g (\infty)$ for some branching geodesic $\g\in \B\G$.
\end{Prop}

\begin{Prf}
We first observe that one implication in the Proposition is
immediate: if $p=\g (\infty)$ for some branching geodesic $\g\in
\B\G$, then $p$ is branching. So let us focus on the converse.

Assume that $p\in \pinf \tX$ is a branching point, and that $p$ is
not a limit point of any branching geodesic.  Let $\gamma$ be a
geodesic ray with $\gamma (\infty)=p$, and observe that for the
geodesic ray $\gamma$ we have either:
\begin{enumerate}
\item $\gamma$ passes through finitely many lifts of chambers, or
\item $\gamma$ passes through infinitely many lifts of chambers.
\end{enumerate}
For each of the two cases above, we need to show that $p$ cannot
be branching. We argue by contradiction.  Assume that
$f:T\rightarrow \pinf \tX$ be an injective mapping of the tripod
into $\pinf \tX$, satisfying $f(*)=p$.

\vskip 5pt

\noindent {\bf Case 1:} Since $\gamma$ passes through finitely
many lifts of chambers, and since the lifts of chambers are
totally geodesic subsets of $\tX$, we have that there exists a
fixed lift $C$ of a chamber with the property that $\gamma (t)\in
C$ for all $t$ sufficiently large. Now pick a basepoint $x$ in the
{\it interior} of $C$, and pick $\epsilon$ small enough so that
the $\epsilon$ metric ball $B_\epsilon (x)$ centered at $x$ is
contained entirely in the interior of $C$. Denote by $lk(x)$ the
boundary of this metric ball, and observe that $lk(x)$ is
homeomorphic to $S^1$. Let $\rho:\pinf \tX\rightarrow lk(x)$ be
the geodesic projection, and consider the composite map $\rho\circ
f:T\rightarrow lk(x)\cong S^1$.

We first note that there is no injective map from $T$ to $S^1$,
hence the composite $\rho\circ f$ must fail to be injective at
some point.  Let $I\subset S^1$ be the subset of points in $lk(x)$
where the map $\rho$ is injective.  Our goal is now to show that
the composite $\rho\circ f$ fails to be injective at some point $z
\in I$.  This immediately implies that $f$ fails to be injective
at the point $\rho ^{-1}(z)$, which would yield the desired
contradiction.

Let us start by observing that $I$ consists of a Cantor set in
$S^1$.  Indeed, if a point $w$ lies in the complement of $I$, then
there exist a pair of geodesic rays $\gamma_1,\gamma_2$ emanating
from $x$, both of which pass through $w\in lk(x)$, but satisfying
$\gamma_1(\infty)\neq \gamma_2(\infty)$.  In particular, the
geodesic rays $\gamma_1, \gamma_2$ must coincide for a period of
time, and subsequently diverge.  This implies that $\gamma_i\cap
\partial C\neq \emptyset$.  Hence the point $w$ lies in the image
$\bar \rho (\partial C)$ of $\partial C$ under the geodesic
retraction map $\bar \rho:\tX- B_\epsilon(x)\rightarrow lk(x)$.
Conversely, given a point $w$ in $\bar \rho (\partial C)$, one can
easily construct a pair of geodesic rays $\gamma_1,\gamma_2$
originating from $x$, passing through $w$, but having
$\gamma_1(\infty)\neq \gamma_2(\infty)$.  This forces the
complement of $I$ to coincide with the set $\bar \rho (\partial
C)$.  But the complement of the set $\bar \rho (\partial C)$ can
naturally be identified with $\pinf C$.  Since $C$ is the
universal cover of a compact negatively curved surface with
non-empty boundary, it is quasi-isometric to a free group $F_k$.
This implies that $I$ is homeomorphic to $\pinf F_k$, which is
known to be a Cantor set.

Now observe that the complement of the set $I\subset lk(x)$
consists of a countable dense union of open intervals.  Let
$I_\partial\subset I$ denote the subset of $I$ consisting of the
boundary points of these intervals.  Note that, by the discussion
above, the set $I_\partial$ coincides with the set $\rho (\pinf
(\partial C))$, and since we are assuming that the point $p\in
\pinf X$ is not the limit point of a branching geodesic, we have
that $(\rho\circ f)(*)=\rho(p)\in I-I_\partial$.

Let $L_i\cong [0,1)$ ($1\leq i\leq 3$) denote the three components
of $T-*$, which we will call the {\it open leaves} of the tripod
$T$.  Since $(\rho\circ f)(*)\in I$, we have that $(\rho\circ
f)(*)\notin (\rho\circ f)(L_i)$ for each $i$.  Let $U$ denote a
small open connected neighborhood of $(\rho\circ f)(*)$ in $lk(x)\cong S^1$, 
and observe that $(\rho\circ f)(*)$ locally separates $U$ into a pair
of open intervals $U_1,U_2$.  If $U$ is chosen small enough, we
must have that a pair of leaves surjects onto one of the $U_j$. We
assume, without loss of generality that $U_1\subset (\rho\circ
f)(L_1)\cap (\rho\circ f)(L_2)$.  But now we note that in the
Cantor set $I$, every point in $I-I_\partial$ can be approximated
on both sides by points in $I_\partial$.  This implies that
$[(\rho\circ f)(L_1)\cap (\rho\circ f)(L_2)]\cap I_\partial \neq
\emptyset$, and as $\rho$ is injective on $I_\partial$, that there
exist a pair of points $q_1\in L_1, q_2\in L_2$ with
$f(q_1)=f(q_2)\in \pinf \tX$.  But this contradicts the fact that
$f$ is injective, concluding the proof for the first case.

\vskip 5pt

\noindent {\bf Case 2:}  For the second case, we assume that $\g$
is a geodesic ray with $\g (\infty)=p$, and which passes through
infinitely many lifts of chambers.  Note that this forces the
geodesic ray $\g$ to intersect infinitely many branching
geodesics.  Let $\{\eta_i\}$ be the collection of branching
geodesics intersected by $\g$, indexed in the order in which their
intersections occur along $\g$.  We now recall two facts:
\begin{enumerate}
\item each $\pinf \eta_i$ separates $\pinf \tX$ (Lemma 3.4 above), 
and \item if $U_i$ denotes the path-connected component of
$\pinf \tX-\pinf \eta_i$ containing $p$, then the collection
$\{U_i\}$ forms an open, path-connected, neighborhood base of $p$
in $\pinf \tX$ (see the proof of Proposition 2.3 in \cite{la1}).
\end{enumerate}

Armed with these two facts, the argument for Case 2 is easy: let
$L_i\cong [0,1)$ again denote the three open leaves of the tripod.
Since $f:T\rightarrow \pinf \tX$ is injective, we have that
$p\notin f(\partial T)$, and hence there exists a small enough
neighborhood $N$ of $p$ with the property that $f(\partial T)
\subset \pinf \tX -\bar N$.  From Fact (2) above, we have that
$U_i\subset N$ for $i$ sufficiently large, and hence that
$f(\partial T)\subset \pinf \tX -\bar N$.  From Fact (1) above, we
have that the corresponding $\pinf \eta_i=\{\eta_i(\pm \infty)\}$
separates $f(\partial T)$ from $p$.  But this implies that, for
each $1\leq i\leq 3$, we have that $f(L_i)\cap\{\eta_i(\pm
\infty)\} \neq \emptyset$.  The pigeon-hole principle forces the
image of a pair of leaves to pass through one of the two points
$\{\eta_i(\pm \infty)\}$.  But this contradicts the fact that $f$
is injective, concluding the argument for Case 2, and completing
the proof of the Proposition.
\end{Prf}

We observe an immediate corollary of the above proposition:

\begin{Cor}
If $f:\tX_1\rightarrow \tX_2$ is a quasi-isometry, and
$f_\infty:\pinf \tX_1\rightarrow \pinf \tX_2$ the induced map on
the boundary at infinity.  Then $f_\infty$ restricts to a
bijection from $\pinf \B\G_1$ to $\pinf \B\G_2$.
\end{Cor}

\begin{Prf}
Note that the induced map $f_\infty: \pinf \tX_1\rightarrow \pinf
\tX_2$ is a homeomorphism. But the previous Proposition
characterizes the subsets $\pinf \B\G_i\subset \pinf \tX_i$ purely
topologically, yielding the corollary.
\end{Prf}

We would now like to focus on the specific situation at hand,
namely we will assume that we are given a pair $X_1,X_2$ of
simple, thick, 2-dimensional hyperbolic P-manifolds, and an
abstract isomorphism $\phi$ from $\Gamma_1:=\pi_1(X_1)$ to
$\Gamma_2:=\pi_1(X_2)$. We observe the following facts that hold
in this setting:
\begin{itemize}
\item the groups $\Gamma_i$ act by homeomorphisms on $\pinf
\tX_i$, \item the isomorphism $\phi$ induces a quasi-isometry
$\bar \phi$ from $\tX_1$ to $\tX_2$, \item the quasi-isometry
$\bar \phi$ induces a homeomorphism $\pinf \bar \phi: \pinf
\tX_1\rightarrow \pinf \tX_2$ which is equivariant with respect to
the $\Gamma_i$ actions on the $\pinf \tX_i$.
\end{itemize}

The previous corollary tells us that $\pinf \bar \phi$ restricts
to a bijection between from the set $\pinf \B\G_1$ to the set
$\pinf \B\G_2$.  Now notice that any branching geodesic $\g \in
\B\G_1$ naturally corresponds to a {\it pair} of points $\{\g (\pm
\infty)\}\subset \pinf \B\G_1$.  We would like to ensure that,
under our homeomorphism $\phi _\infty$, the pair $\{\g (\pm
\infty)\}$ maps to a pair $\{\g^\prime (\pm \infty)\}$ for some
branching geodesic $\g^\prime \in \B\G_2$.  In order to achieve
this, our next step is to characterize the endpoints of branching
geodesics in a purely topological manner. This is the content of
our:

\begin{Prop}
Let $\{p,q\}\subset \pinf \B\G\subset \pinf \tX$ be an arbitrary
pair of distinct points.  Then we have that:
\begin{enumerate}
\item if there exists a $\g \in \B\G$ with the property that $\{\g
(\pm \infty)\}=\{p,q\}$, then $\{p,q\}$ separates $\pinf \tX$ into
$\geq 3$ components. \item if there exists a geodesic $\g$
contained in the interior of a single lift of a chamber, with the
property that $\{\g (\pm \infty)\}=\{p,q\}$, then $\{p,q\}$
separates $\pinf \tX$ into exactly $2$ components. \item in all
other cases, $\{p,q\}$ does not separate $\pinf \tX$.
\end{enumerate}
\end{Prop}

\begin{Prf}
Statement (1) is an immediate consequence of Lemma 3.4 and the
thickness hypothesis.

To see statement (2), one starts with a $\g\notin \B\G$, and
satisfying $\g \subset Int(C)$ with $C\in \C$.  Note that this
implies that $\g$ separates $C$ into precisely two open components,
denoted $Z_1$ and $Z_2$. Furthermore, the closure of each
component is a closed, totally geodesic subset of $\tX$.  Now for
$i=1,2$, define the sets $(Z_i)_j$ ($j\geq 1$) inductively by:
\begin{itemize}
\item the initial condition $(Z_i)_1=Z_i$, and \item $(Z_i)_{j+1}$
is the union of $(Z_i)_j$ along with all lifts of chambers which
are incident to $(Z_i)_j$.
\end{itemize}
We observe that we have proper inclusions $(Z_i)_j\subset
(Z_i)_{j+1}$, and that each of the subsets $(Z_i)_j$ is totally
geodesic and path-connected.

Now form the sets $Y_i:=\cup_{j\in \mN} (Z_i)_j$, and observe that
each $Y_i$ is a path-connected, totally geodesic subspace of $\tX$
(as the latter properties are preserved under increasing unions).
Furthermore each of the sets $Y_i$ has the property that $\pinf
Y_i-\{\g(\pm\infty)\}$ is path-connected.  Indeed, given a pair of
points in $\pinf Y_i$, one can consider the geodesic $\eta$
corresponding to the pair of points.  It is easy to see that there
is an isometrically embedded, totally geodesic ``half-$\mH^2$''
$H\subset Y_i$ with boundary the given geodesic $\eta$.  This
implies that there exists $\pinf H\cong \mI\subset \pinf
Y_i-\{\g(\pm\infty)\}$ whose endpoints correspond precisely to
$\{\eta (\pm\infty)\}$.

Finally, observe that $\tX-\g=Y_1\coprod Y_2$, and that the
closure of $Y_i$ is precisely $Y_i\cup \g$.  Hence we have that
$\pinf \tX=\pinf Y_1 \cup_{\{\g(\pm\infty)\}}\pinf Y_2$,
expressing $\pinf \tX$ as a union of a pair of closed sets (as
each $Y_i$ is totally geodesic) whose intersection is precisely
$\{\g(\pm\infty)\}$, and with the property that each $\pinf
Y_i-\{\g(\pm\infty)\}$ is path-connected.  This immediately
implies statement (2) of the proposition.

So we are now left with showing statement (3).  In order to do this,
we first make two observations concerning branching geodesics.  Note
that if $\rho \in \B\G$, then we have that $\pinf \tX -\{\rho (\pm \infty)\}$
splits into $k\geq 3$ path-connected components $U_1,\ldots ,U_k$.
We now observe:

\vskip 5pt

\noindent {\bf Fact 1:}  The closure of each $U_i$ is $\bar U_i=U_i\cup 
\{\rho (\pm \infty)\}$, and is path-connected.

\vskip 5pt

\noindent {\bf Fact 2:}  For every pair of distinct points $x,y\in Y_i$ there is 
a path $\eta_x\subset \bar U_i -y$ joining $x$ to one of the points $\{\rho (\pm 
\infty)\}$.

\vskip 5pt

Now assuming these two facts, we proceed with the proof of statement (3).
For the points $\{p,q\}$ satisfying the hypotheses of statement (3), we
have that the geodesic $\g$ satisfying $\g (\pm \infty)=\{p,q\}$ must 
intersect a branching geodesic $\rho \in \B\G$.  Upto re-indexing, 
we may assume that $p \in U_1$, and $q\in U_2$.  

Now pick an arbitrary pair of points $\{x,y\} \subset \pinf \tX -\{p,q\}$, 
and consider the subsets $U_i,U_j$ satisfying $x\in U_i$, $y\in U_j$.  If $i\geq 3$,
then from {\bf Fact 1}, there exists a path in $\bar U_i$ joining $x$ to $\rho (\infty)$.
On the other hand, if $i=1,2$, then from {\bf Fact 2}, there is a path joining $x$ to
one of the points $\rho (\pm \infty)$, avoiding the point $p$ (if $i=1$) or $q$ (if 
$i=2$).  In either case, denote this path by $\eta_x$.  Now applying the same 
reasoning to $y$, we find a path $\eta_y$ joining $y$ to one of the points $\rho(\pm
\infty)$, and avoiding the pair $\{p,q\}$.  

If the endpoints of the paths $\eta_x,\eta_y$ coincide, concatenation gives us 
a path connecting $x$ to $y$.  Otherwise, from Fact 1, we note that there is a
path $\eta_{\rho}$ in $Y_3$ joining $\rho(\infty)$ to $\rho (-\infty)$.  
Concatenating the three paths $\eta_x$, $\eta_\rho$, and $\eta_y$ yields a path
joining $x$ to $y$ in $\pinf \tX-\{p,q\}$.  Since this holds for arbitrary $x,y\in
\pinf \tX-\{p,q\}$, we conclude that $\pinf \tX-\{p,q\}$ is path-connected.  So
to complete the proof of the Proposition, we are left with verifying {\bf Fact 1} 
and {\bf Fact 2}.

To see {\bf Fact 1}, we first note that the complement of $U_i\cup \{\rho 
(\pm \infty)$ consists of the union $\coprod_{j\neq i} U_j$.  Since all the $U_j$
are open, this implies that the closure of $U_i$ is contained in the set 
$U_i\cup \{\rho (\pm \infty)$.  To see the converse, we observe that we can 
construct, as in the argument for statement (2), totally geodesic subspaces $Y_i$
with the property that $\pinf Y_i=U_i\cup \{\rho (\pm \infty)\}$ and with 
$\partial Y_i=\rho$.  But within
the $Y_i$, it is easy to see that there exist totally geodesic embedded `half' 
$\mH^2$'s whose boundary is precisely $\rho$.  At the level of the boundary at 
infinity, this yields an embedded interval in $U_i\cup \{\rho (\pm \infty)\}$
with endpoints precisely $\{\rho (\pm \infty)\}$.  This immediately implies
that $\{\rho (\pm \infty)\}$ lies in the closure of the $U_i$, completing the
proof of {\bf Fact 1}.

To see {\bf Fact 2}, we note that given $p\in U_i$, the embedded interval 
$f:\mI \hookrightarrow U_i\cup \{\rho (\pm \infty)\}$ mentioned in the previous
paragraph can be chosen to satisfy $f(0)=\rho(+\infty)$, $f(1)=\rho(-\infty)$,
and $f(1/2)=p$.  Now note that if $q\notin f(\mI)$, we are done.  On the other
hand, if $q\in f(\mI)$, then the hypothesis that $p\neq q$ ensures that $q=f(r)$
where either $0<r<1/2$ or $1/2<r<1$.  In both cases we can use $f$ restricted to
a suitable subinterval of $\mI$ to get the desired path.  This completes the
proof of {\bf Fact 2}, and hence, of Proposition 3.2.

\end{Prf}

Since separation properties are purely topological, we obtain the
immediate corollary:

\begin{Cor}
Every quasi-isometry $f:\tX_1\rightarrow \tX_2$ naturally induces
a bijective correspondance between $\B\G_1$ and $\B\G_2$.
\end{Cor}

Now observe that for the universal cover $\tX$ of a simple, thick,
2-dimensional hyperbolic P-manifold, we can naturally define an
adjacency relation on the set $\B\G$.  We say that a pair of
elements $\gamma_1, \gamma_2$ of $\B\G$ are {\it adjacent},
denoted by $\gamma_1 \sim \gamma_2$ provided there exists a
geodesic joining a point in $\gamma_1$ to a point in $\gamma _2$,
and {\it lying entirely within a single chamber}. Note that the
above relation is symmetric, but {\it not} transitive.  The next
step is to establish that a quasi-isometry preserves the adjacency
relation.

\begin{Prop}
If $f:\tX_1\rightarrow \tX_2$ is a quasi-isometry, then for any
pair $\gamma_1,\gamma_2 \in \B\G_1$, we have:
$$\gamma_1 \sim \gamma_2 \Longleftrightarrow f_*(\gamma_1) \sim f_*(\gamma_2),$$
where $f_*(\g)$ in $\B\G_2$ is the branching geodesic
bijectively associated with $\g\in \B\G_1$.
\end{Prop}

\begin{Prf}
This follows immediately from Proposition 3.2.  Indeed, from the
definition of the relation $\sim$, we see that $\g _1\sim \g _2$
if and only if the geodesics form a pair of distinct boundary geodesics of
a single chamber $C$.  Given a 4-tuple of distinct points $\{x_-, x_+,y_-,y_+\}
\subset \pinf \tX$, Proposition 3.2 (parts (1) and (2)) tells us there is
a pair $\g _1, \g_2$ satisfying $\g_1(\pm \infty)=x_\pm$, $\g_2(\pm \infty)
=y_\pm$ if and only if we have:
\begin{itemize}
\item the two pairs of points $\{x_\pm\}$, $\{y_\pm\}$ each separate $\pinf \tX$
into $\geq 3$ components,
\item each of the four pairs of points $\{x_+,y_+\}$, $\{x_+,y_-\}$, 
$\{x_-,y_+\}$, $\{x_-,y_-\}$, separate $\pinf \tX$ into precisely two components.
\end{itemize}

Since the quasi-isometry $f$ induces a homeomorphism $f_\infty$ between the
two boundaries at infinity $\pinf \tX_1$ and $\pinf \tX_2$, the above topological
characterization of endpoints of adjacent branching geodesics immediately yield
the Proposition. 
\end{Prf}

Next observe that the adjacency relation can be used to keep track
of the chambers. This is the content of the following:

\begin{Lem}
Elements of $\C$ correspond bijectively to maximal subsets of
$\B\G$ on which the adjacency relation is transitive.
\end{Lem}

\begin{Prf}
Let $C\in \C$ be a chamber, and associate to it the collection of 
$\g \in \B\G$ which arise as the boundary components of $C$; denote this set
by $B_C$.  It is immediate from the definition of the relation $\sim$ that 
the adjacency relation is transitive on $B_C$.

Conversely, let $B\subset \B\G$ be a subset on which the adjacency relation
is transitive.  We claim that there is a $C\in \C$ satisfying $B\subset B_C$.
To see this, pick $\g _1, \g _2, \g _3\in B$, and observe that the condition
$\g_1\sim \g_2$ implies that the two branching geodesics are boundary components
of a fixed chamber $C_{12}$.  Similarly, $\g_2\sim \g_3$ implies that they are
both boundary components of a chamber $C_{23}$.  Now note that if $C_{12}\neq 
C_{23}$, then they form two distinct chambers both incident to $\g _2$.  From
Lemma 3.2, this forces $\g_1$ and $\g_3$ to lie in distinct connected components
of $\tX -\g_2$.  Hence, if $\eta$ is an arbitrary geodesic segment joining a 
point on $\g_1$ to a point on $\g_3$, we have that $\eta \cap \g_2\neq \emptyset$.
Since $\eta$ is assumed to be geodesic, by restricting we can view $\eta$ as 
a concatenation of a geodesic joining a point in $\g_1$
to a point in $\g_2$, together with a geodesic joining a point in $\g_2$ to a 
point in $\g_3$.  Since $\g_1,\g_2$ are distinct boundary components of $C_{12}$,
the first geodesic segment must intersect the interior of $C_{12}$ non-trivially.  
Likewise, the second geodesic segment must intersect the interior of $C_{23}$ 
non-trivially.  But this contradicts the assumption that $\g_1\sim \g_3$.

Since the adjacency relation is transitive on the subsets $B_C$ , and since 
every subset on which the adjacency relation is transitive is contained in one of
the $B_C$, we conclude that the latter are precisely the maximal subsets on 
which $\sim$ is transitive, concluding the proof of the Lemma.

\end{Prf}

By combining the previous Lemma with the previous proposition, we
immediately obtain the:

\begin{Cor}
If $f:\tX_1 \rightarrow \tX_2$ is a quasi-isometry, then $f$
induces a bijection $f_*$ from $\C_1$ to $\C_2$.  Furthermore, if $\gamma \in
\B\G_1, C\in \C_1$ satisfy $\gamma \subset C$, then $f_*(\gamma) \subset
f_*(C)$.
\end{Cor}

We now know that a quasi-isometry between the universal covers of
a pair of simple, thick, 2-dimensional hyperbolic P-manifolds
induces a bijection between the lifts of chambers.  Now recall
that in the situation we are interested in, the quasi-isometry
$\bar \phi_\infty$ from $\tX_1$ to $\tX_2$ has the additional
property that it is $(\Gamma_1,\Gamma_2)$-equivariant.  In
particular, each lift of a chamber $C\in \C_i$ has a stabilizer
under the action of $\Gamma_i$ on $\tX_i$.

Our next step is to identify the stabilizers of the various $C\in
\C_i$ from the boundary at infinity.  This is made precise in the
following:

\begin{Prop}
Consider the natural action of $\Gamma$ on $\tX$ (where
$\Gamma=\pi_1(\tX)$), and the corresponding induced action on
$\pinf \tX$. Then for every $C \in \C$ we have that the stabilizer
of $C$ under the $\pi_1(C)$-action on $\tX$ coincides with the 
stabilizer of $\pinf C$ under the induced $\pi_1(C)$-action on
$\pinf \tX$.
\end{Prop}

The argument for this Proposition can be found in \cite{la1} (see
the Assertion on pg. 212).  
Since an abstract isomorphism between $\Gamma_1$ and $\Gamma_2$
yields a $(\Gamma_1,\Gamma_2)$-equivariant homeomorphism between
$\pinf \tX_1$ and $\pinf \tX_2$, this immediately yields the:

\begin{Cor}
The bijection $\bar \phi_*:\C_1 \rightarrow \C_2$ induced by the
isomorphism $\phi$ has the property that, for every $C\in \C_1$,
one has $Stab_{\Gamma_1}(C)\cong Stab_{\Gamma_2}(\bar \phi_*(C))$.
\end{Cor}

Armed with this information, it is now easy to complete the proof
of the topological rigidity Theorem 1.2.  We first note that, by
the $(\Gamma_1,\Gamma_2)$-equivariance of the homeomorphism from 
$\pinf \tX_1$ to $\pinf \tX_2$, and in view of Proposition 3.3, 
we conclude that there are an equal number of orbits of branching
geodesics in $\tX_1$ and $\tX_2$.  Since each such orbit corresponds
to precisely one circle in the branching locus of the respective
$X_i$, we conclude that there is a bijective correspondance between
the components of the branching locus of $X_1$ and $X_2$.

Now given any chamber in $X_1$, one can consider the family of 
connected lifts of the chamber.  These form a single orbit under
the $\Gamma_1$-action on $\C_1$.  Each lift in this orbit, by 
Corollary 3.4 maps to the lift of a corresponding chamber in $\C_2$.
But $(\Gamma_1,\Gamma_2)$-equivariance ensures that the lifts of
chambers one gets in $\C_2$ lie in a single $\Gamma_2$-orbit under
the corresponding $\Gamma_2$-action on $\C_2$.  Finally, Corollary
3.5 implies that the original chamber in $X_1$, and the corresponding
chamber in $X_2$ have the same fundamental group.  

Next we note that both the number of boundary components the chamber
has can be entirely determined by the number of distinct orbits of 
of points $\{\gamma (\pm \infty)\}$ lying in $\pinf C$, where $C\in \C_i$
is a lift of the chamber.  Again, by $(\Gamma_1,\Gamma_2)$-equivariance
of the homeomorphism, we obtain that pairs of corresponding chambers in
$X_1$ and $X_2$ have exactly the same number of boundary components.  Now
observe that a surface with boundary is uniquely determined by it's fundamental
group and the number of boundary components it has.  This tells us that
the correspondance between chambers in $X_1$ and $X_2$ preserves the
homeomorphism type of the chambers.  To conclude, we note that the dynamics
on the boundaries at infinity also allow us to keep track of how each chamber
is attached to the branching strata.  Putting all this together, we obtain
a homeomorphism from $X_1$ to $X_2$.  It is immediate by construction that
this homeomorphism induces the original isomorphism on the level of the 
fundamental groups, completing the proof of Theorem 1.2.

\section{Concluding remarks.}

Our main theorem states that, within a certain class of diagrams of 
groups, each group that appears as a direct limit has a unique 
representative.  An interesting question is the following:

\vskip 5pt

\noindent {\bf Question:}  Which classes of diagrams of groups have
the property that any group that occurs as a direct limit arises as
the limit of a unique diagram?

\vskip 5pt

Forester \cite{F} has given criterions under which a Bass-Serre 
splitting of a group is unique (see also Guirardel \cite{G}).  We note
that the Bass-Serre trees naturally associated to geometric amalgamations
of free groups {\it do not} satisfy the hypotheses in Forester's work.

Another interesting aspect of the groups we are considering lies in 
the fact that these groups are essentially combinatorially determined.
Indeed, in order to recognize the isomorphism type of these groups, 
it is sufficient (by the main theorem) to keep track of:
\begin{itemize}
\item the ranks of the free groups arising as fundamental groups of chambers, 
\item the number of boundary components of each chamber
\item how the chambers get glued
\end{itemize}
Since this information consists of a finite amount of data, these
groups form a class of $\delta$-hyperbolic groups (or $CAT(-1)$ groups) 
for which one can look at various decision type problems.  We can ask:

\vskip 5pt

\noindent {\bf Question:} When are the direct limits of a pair of 
geometric amalgamations of free groups quasi-isometric?  When are they
bi-Lipschitz equivalent?

\end{document}